\documentclass{article}

\usepackage{zharticle}
\usepackage{graphicx,amsmath,amsfonts,verbatim}
\usepackage{enumitem}
\numberwithin{equation}{section}
\allowdisplaybreaks%

\newcommand{\ones}{{\bf 1{}}}  % vector with all components one
\newcommand{\reals}{{\bf R{}}}  % real numbers
\newcommand{\prob}{\mathop{\bf prob{}}}  % probability
\newcommand{\card}{\mathop{\bf card}}  % cardinality
\newcommand{\tr}{\mathop{\bf tr}}  % trace
\newcommand{\argmax}{\mathop{\rm argmax}}  % argmax
\newcommand{\norm}[1]{{\lVert #1 \rVert}}  % norm
\newcommand{\cvxset}{{\cal C{}}}  % convex set
\newcommand{\idxset}{{\cal F{}}}  % index set
\newcommand{\Dkl}{\mathop{D_{\rm kl}}}  % KL-divergence

\newcommand{\cf}{{\rmfamily\itshape cf.}}
\newcommand{\eg}{{\rmfamily\itshape e.g.}}
\newcommand{\ie}{{\rmfamily\itshape i.e.}}
\newcommand{\etc}{{\rmfamily\itshape etc.}}
\newcommand{\etal}{{\rmfamily\itshape et al.}}

\usepackage{url}
\usepackage[colorlinks]{hyperref}

\usepackage{lipsum}

\usepackage{listings}
\usepackage{xcolor}

\definecolor{codegreen}{rgb}{0,0.6,0}
\definecolor{codegray}{rgb}{0.5,0.5,0.5}
\definecolor{codepurple}{rgb}{0.58,0,0.82}
\definecolor{backcolour}{rgb}{0.95,0.95,0.92}

\lstdefinestyle{zhstyle}{
    commentstyle=\color{codegreen},
    keywordstyle=\color{magenta},
    numberstyle=\tiny\color{codegray},
    stringstyle=\color{codepurple},
    basicstyle=\ttfamily\small,
    breakatwhitespace=false,
    breaklines=true,
    captionpos=b,
    keepspaces=true,
    numbers=left,
    numbersep=5pt,
    showspaces=false,
    showstringspaces=false,
    showtabs=false,
    tabsize=2,
    frame=topbottom,
    framerule=0.5pt
}
\lstdefinelanguage{zhpython}{%
    sensitive=true,%
    morekeywords={and, as, assert, break, class, continue, def, del, elif,%
        else, except, exec, finally, for, from, global, if, import, in, is,%
        lambda, not, or, pass, print, raise, return, try, while, with, yield},%
    morecomment=[l]\#,%
    morestring=[s]{'''}{'''},%
    morestring=[s]{"""}{"""},%
    morestring=[b]',%
    morestring=[b]"%
}

\title{Multi-convex Programming for Discrete\\Latent Factor Models Prototyping}
\author{Hao Zhu}
\author{Shengchao Yan}
\author{Jasper Hoffmann}
\author{Joschka Boedecker}
\affil{Department of Computer Science and IMBIT//BrainLinks-BrainTools\\University of Freiburg}

\begin{document}
\maketitle

\begin{abstract}
    Discrete latent factor models (DLFMs) are widely used in various domains such as machine learning, economics, neuroscience, psychology, \etc\
    Currently, fitting a DLFM to some dataset relies on a customized solver for individual models, which requires lots of effort to implement and is limited to the targeted specific instance of DLFMs.
    In this paper, we propose a generic framework based on \texttt{CVXPY}, which allows users to specify and solve the fitting problem of a wide range of DLFMs, including both regression and classification models, within a very short script.
    Our framework is flexible and inherently supports the integration of regularization terms and constraints on the DLFM parameters and latent factors, such that the users can easily prototype the DLFM structure according to their dataset and application scenario.
    We introduce our open-source Python implementation and illustrate the framework in several examples.
\end{abstract}

\newpage
\tableofcontents
\newpage

\section{Introduction}
Discrete latent factor models (DLFMs) represent a class of models in which some variables are always hidden.
In this paper, we consider DLFMs expressed as,
\begin{equation}\label{eq:dlfm_std}
    \begin{array}{rcl}
        z & \sim & \prob(z)\\
        y & \sim & \prob(y \mid x, z, \theta),
    \end{array}
\end{equation}
where the \emph{latent factor} $z \in \{1, \ldots, K\}$ is a discrete random variable, $x$ and $y$ are the \emph{feature} and \emph{observation}, respectively, and $\theta$ denotes the \emph{model parameters}.
For example, a \emph{mixture of linear regressions} in standard form (\ref{eq:dlfm_std}) is
\begin{equation*}
    \begin{array}{rcl}
        z & \sim & {\rm Cat}(p)\\
        y & \sim & {\cal N}(x^T \theta_z, \sigma^2)
    \end{array}
\end{equation*}
with latent factor $z \in \{1, \ldots, K\}$, feature $x \in \reals^n$, observation $y \in \reals$, and parameters $\theta_z \in \reals^n$ for all $z = 1, \ldots, K$, where ${\rm Cat}(p)$ denotes the categorical distribution with category probability vector $p \in \reals^K_+$ satisfying $\ones^T p = 1$, and ${\cal N}(x^T \theta_z, \sigma^2)$ is the Gaussian distribution with mean $x^T \theta_z$ and variance $\sigma^2$.
DLFMs appear in domains such as machine learning, economics, signal processing, and cognitive science, with a vast range of applications~\cite[chapter 28]{murphy2023probabilistic}.
Specifically, in neuroscience and psychology, DLFMs have provided interpretable characterizations of both neural population activities~\cite{jha2021factor} and subject behavior~\cite{ashwood2022mice,zhu2024multiintention}.
Note that the references listed here are by no means thorough; interested readers may refer to Jha \etal~\cite{jha2024active} for a detailed review.
In \S\ref{sec:example}, we also provide some specific examples applied in these domains.

\paragraph{DLFM fitting problems.}
We consider the problem of fitting a DLFM model given the set of $m$ observations $\{y_1, \ldots, y_m\}$ and the corresponding features $\{x_1, \ldots, x_m\}$.
Informally, the DLFM fitting problem consists in finding the parameters $\theta_z$ for all $z = 1, \ldots, K$, and the latent factors $z_i$ for all $i = 1, \ldots, m$, such that some metric indicating the model fitting error (\eg, the negative log-likelihood of the data) is minimized.
(We provide a formal definition of the DLFM fitting problem in \S\ref{sec:dlfm_fit_prob_std}.)
The DLFM fitting problem is in general very hard to solve (actually, NP-hard, as we will show in \S\ref{sec:dlfm_fit_prob_std}).
Since the latent factors $\{z_1, \ldots, z_m\}$ dominating the generation of the observations $\{y_1, \ldots, y_m\}$ under features $\{x_1, \ldots, x_m\}$ were not observed, in most cases, fitting a DLFM requires solving a nonconvex combinatorial optimization problem where the number of variables increases linearly with the number of samples in the dataset.

\paragraph{Current methods for fitting DLFMs.}
Various approaches have been introduced to fit different classes of DLFMs.
The most commonly used method for fitting DLFMs is the \emph{expectation-maximization} (EM) algorithm~\cite{dempster1977maximum}.
EM algorithm is an iterative method with each iteration alternating between performing an expectation (E) step, which creates a function for the expectation of the log-likelihood evaluated using the current estimate for the parameters, and a maximization (M) step, which computes parameters maximizing the expected log-likelihood found on the E-step.
The estimated parameters are then used to determine the distribution of the latent factors in the next E-step.
The major limitation with EM is that only convergence to a local minimum is guaranteed.
Thus, it is common to apply repeated EM iterations with random initializations and select the best-performing estimation, which turns out to work quite well in practice~\cite{ashwood2022mice,zhu2024multiintention}.
Another class of methods for fitting DLFMs is the \emph{Monte Carlo} methods~\cite[chapters 11--13]{murphy2023probabilistic}, which instead of directly solving the fitting problem, but try to estimate the posterior distribution of the unknown parameters given the problem data and some predefined prior, by generating a large number of instances from the desired distribution.
Monte Carlo methods avoid the local optimality limitation of EM methods and are hence more accurate than EM, but normally require a very long time to generate enough instances for an accurate estimation of the targeted posterior distribution.
A compromise between EM and Monte Carlo is the \emph{variational inference} methods~\cite[chapter 10]{murphy2023probabilistic}, which also try to (approximately) estimate a posterior distribution of the unknown parameters as in Monte Carlo methods, but via solving an alternative optimization problem.
A more detailed discussion about Monte Carlo and variational inference methods for fitting DLFMs can be found in Jha \etal~\cite{jha2024active}.

\paragraph{Limitations.}
The aforementioned methods only provide an outline for designing algorithms, specifically, each for a different DLFM.
To fully implement an algorithm instance that can be directly applied to the given data to fit a DLFM, one would need expert knowledge in lots of domains, such as probability, statistics, linear algebra, optimization, \etc, as well as vast minor but essential coding tricks.
Although stable implementation of classic DLFMs (\eg, hidden Markov models) is provided by open-source packages such as \texttt{scikit-learn}~\cite{pedregosa2011scikitlearn}, the more recently developed DLFMs are published only (and hopefully) accompanied with source code closely dependent on the respective practical problems.
In the former case, it would be challenging to integrate even minor adaptations to the provided DLFMs, such as constraints on model parameters.
In the latter case, understanding the open-source implementation for the DLFM would already require a lot of effort, and at least medium level knowledge in the related domains.
For the users, this blocks them from applying the powerful class of DLFMs, especially the more recently developed ones, to practical problems, since before launching the machinery designed for very fast and/or accurate DLFM fitting, the users might first just want to try different models and select the most suitable one for their specific application scenario, and sometimes incorporate minor adaptations.
For the developers, upon proposing a new DLFM or some modifications to the existing DLFMs, lots of effort needs to be devoted to modifying or redesigning the algorithm for solving a different fitting problem, before even knowing whether the new idea will work or not.
The current status strongly limits both further development and broader application of DLFMs.

\paragraph{Focus of the paper and contributions.}
The focus of this paper is \emph{not} on methods for fitting specific DLFMs, but on a \emph{framework} for DLFMs prototyping, which aims at avoiding the limitations of fitting different DLFMs to some given dataset via respective custom solvers.
Our framework allows users to specify and solve the fitting problem of a wide range of DLFMs easily in a high level, human readable language, as well as introduce custom constraints and regularization terms to the model parameters and latent factors, based on their respective application scenarios.
Our framework supports a wide range of regression and classification DLFMs, whose loss functions and parameter constraints are representable via \emph{disciplined convex programming} (DCP), even though the fitting problem itself is not convex.
Our work can be considered as an extension of the widely used \emph{domain specific language} (DSL), \texttt{CVXPY}~\cite{diamond2016cvxpy}, for specifying, canonicalizing, and solving convex optimization problems, to deal with nonconvex DLFM fitting problems.
Our implementation is fully open-source and can be found at
\begin{quote}
    \url{https://github.com/nrgrp/dlfm}.
\end{quote}

\paragraph{Paper structure.}
We introduce some background knowledge about convex optimization and multi-convex programming in \S\ref{sec:pre}, which is referred to in the sequel.
The standard form of the DLFM fitting problem is introduced in \S\ref{sec:dlfm_fit_prob}.
In this section, we also discuss the properties of the DLFM fitting problem, as well as some examples.
We propose a generic method for (approximately) solving the DLFM fitting problem in \S\ref{sec:heuristic}, by introducing convex relaxations and problem transformation.
In \S\ref{sec:impl}, we provide the implementation of our framework for DLFM prototyping.
Finally, in \S\ref{sec:example}, we list a number of numerical examples.
The objective of showing these examples is \emph{not} on competitive results in terms of solving time, quality, \etc, but rather to show the simplicity of specifying and fitting different DLFMs via our framework, along with results that are at least comparable to those obtained via specialized solvers for the individual problems.

\section{Preliminaries}\label{sec:pre}
\subsection{Disciplined convex programming}\label{sec:dcp}
Disciplined convex programming is a framework for modeling convex optimization problems introduced by Grant \etal~\cite{grant2004disciplined,grant2006disciplined}.
DCP imposes a ruleset, by following which the specified mathematical optimization problem can be easily verified as convex, and then canonicalized to a cone program which is eventually processed by some conic solver.
The conforming problems from DCP are called \emph{disciplined convex programs}.

The DCP ruleset restricts the set of functions that can appear in a problem and the way functions can be composed.
Functions that appear in a disciplined convex program are restricted to those that are some composition of a set of atomic functions with known curvature and graph implementation, or can be represented as partial optimization over a cone program~\cite{nesterov1992conic,grant2008graph}.
Suppose we are given some function $f = h(g_1(x), \ldots, g_p(x))$ where $h \colon \reals^p \to \reals$ is a convex function and $g_1, \ldots, g_p \colon \reals^n \to \reals$, and let $\tilde{h} \colon \reals^p \to \reals \cup \{\infty\}$ be the extended-value extension of $h$~\cite[\S3.1]{boyd2004convex}.
The function $f$ is convex if and only if one of the following conditions is satisfied for all $i = 1, \ldots, p$:
\begin{itemize}
    \item $g_i$ is convex and $\tilde{h}$ is nondecreasing in the $i$th argument;
    \item $g_i$ is concave and $\tilde{h}$ is nonincreasing in the $i$th argument;
    \item $g_i$ is affine.
\end{itemize}
(The composition rule for concave functions is analogous; see~\cite[\S3.2]{boyd2004convex}.)

A mathematical optimization program has the general form
\begin{equation}\label{prob:dcp_general}
    \begin{array}{ll}
        \mbox{minimize} & f(x)\\
        \mbox{subject to} & g_i(x) \sim h_i(x),\quad i = 1, \ldots, m,
    \end{array}
\end{equation}
where $x \in \reals^n$ is the variable and the relational symbol $\sim$ denotes one of the relational operators $=$, $\leq$, or $\geq$.
Problem (\ref{prob:dcp_general}) is a (DCP-supported) convex optimization problem if the functions $f$, $g_i$, and $h_i$ are expressions with curvature verified by the DCP ruleset and the following curvature restrictions on these expressions are satisfied:
\begin{itemize}
    \item $f$ is convex (concave if (\ref{prob:dcp_general}) is a maximization problem);
    \item when the relational operator is $=$, $g_i$ and $h_i$ are both affine;
    \item when the relational operator is $\leq$, $g_i$ is convex, and $h_i$ is concave;
    \item when the relational operator is $\geq$, $g_i$ is concave, and $h_i$ is convex.
\end{itemize}
Note that an affine expression (function) is both convex and concave, and hence matches either curvature requirement.

Modeling systems for convex optimization problems based on the DCP ruleset have been well developed over the years.
There are multiple DSLs designed for different programming languages, such as \texttt{YALMIP}~\cite{lofberg2004yalmip} and \texttt{CVX}~\cite{grant2014cvx} for MATLAB, \texttt{CVXPY}~\cite{diamond2016cvxpy,agrawal2018rewriting} for Python, and \texttt{Convex.jl}~\cite{udell2014convex} for Julia.
Users are able to specify and solve a DSL verified convex optimization problem in a high level, human readable language, close to the corresponding mathematical expressions.

Our framework introduced in this paper is based on the DSL \texttt{CVXPY}, and uses the aforementioned DCP ruleset to verify whether an instance of the DLFM fitting problem is supported.

\subsection{Multi-convex problems}
We follow the notation from Shen \etal~\cite{shen2017disciplined} for \emph{multi-convex} optimization problems.

\paragraph{Fix variables in a function.}
Consider a function $f \colon \reals^n \to \reals$, and a partition of its variable $x \in \reals^n$ into blocks:
\begin{equation}\label{eq:var_partition}
    x = (x_1, \ldots, x_N),\quad x_i \in \reals^{n_i},\quad \sum_{i = 1}^N n_i = n.
\end{equation}
Let $\idxset \subseteq \{1, \ldots, N\}$ denote an index set, such that $x_\idxset \in \{x_1, \ldots, x_N\}$, and let $\idxset^c = \{1, \ldots, N\} \setminus \idxset$ be the complement of $\idxset$.
By saying \emph{the function $f$ with $\idxset$ fixed at point $\tilde{x} \in \reals^n$}\footnote{Formally, this should be expressed as: \emph{the function $f$ with the variables $x_i$ for $i \in \idxset$ fixed at point $\tilde{x} \in \reals^n$}, but we will follow the convention from Shen \etal~\cite{shen2017disciplined} to use the less formal expression in the text through the paper for simplicity.}, we refer to the function $\tilde{f}$ with variables $x_i$ for $i \in \idxset^c$ and for $i \in \idxset$, $x_i = \tilde{x}_i$.
For example, the function $f(x_1, x_2) = x_1 x_2$ with $\idxset = \{2\}$ fixed at point $\tilde{x} = (3, 5)$ is the function $\tilde{f}(x_1) = 5 x_1$.
We sometimes omit the expression `at point $\tilde{x} \in \reals^n$' to refer to the general case where the function $f$ is fixed at some point.

\paragraph{Multi-convex problems.}
Given an optimization problem with general structure (\ref{prob:dcp_general}) of variable $x \in \reals^n$ partitioned into blocks as in (\ref{eq:var_partition}).
Let $\idxset \subseteq \{1, \ldots, N\}$ be some index set.
We refer to \emph{the problem (\ref{prob:dcp_general}) with $\idxset$ fixed} as the problem
\begin{equation}\label{prob:dcp_general_fix}
    \begin{array}{ll}
        \mbox{minimize} & \tilde{f}(x_{\idxset^c})\\
        \mbox{subject to} & \tilde{g}_i(x_{\idxset^c}) \sim \tilde{h}_i(x_{\idxset^c}),\quad i = 1, \ldots, m
    \end{array}
\end{equation}
with variable $x_{\idxset^c}$, where the functions $\tilde{f}$, $\tilde{g}_i$, and $\tilde{h}_i$, $i = 1, \ldots, m$ are the functions $f$, $g_i$, and $h_i$ with $\idxset$ fixed.
We say \emph{the problem (\ref{prob:dcp_general}) is convex with $\idxset$ fixed}, if the fixed problem (\ref{prob:dcp_general_fix}) is convex, \ie, satisfies the DCP ruleset in \S\ref{sec:dcp}.
We say the problem (\ref{prob:dcp_general}) is \emph{multi-convex}, if there are sets $\idxset_1, \ldots, \idxset_p$, such that for all $i = 1, \ldots, p$, the problem (\ref{prob:dcp_general}) is convex with $\idxset_i$ fixed, and $\bigcap_{i = 1}^p \idxset_i = \emptyset$.
This definition includes convex problems as special cases with $p = 0$, \ie, $\idxset = \emptyset$.
A \emph{biconvex} problem is multi-convex with $p = 2$.
For example, the problem given by
\begin{equation*}
    \begin{array}{ll}
        \mbox{minimize} & |x_1 x_2|\\
        \mbox{subject to} & x_1 + x_2 \geq 1
    \end{array}
\end{equation*}
with variable $x \in \reals^2$ is \emph{not} convex, but is biconvex with index sets $\idxset_1 = \{1\}$ and $\idxset_2 = \{2\}$.

\section{The DLFM fitting problem}\label{sec:dlfm_fit_prob}
Fitting a DLFM to some given dataset is the primary objective and the first step for DLFM users.
In this section, we provide a formal definition to the fitting problem of a wide range of DLFMs, followed by some analysis about the basic properties of such a problem, interpretations, and some specific examples that are widely used.

For simplicity of notation (and with slight abuse of notation), we will denote the latent factor as $z \in \{e_1, \ldots, e_K\} \subseteq \reals^K$ in the sequel, where $e_i$ are the $i$th standard basis vector in $\reals^K$.
Such a vector form notation of the latent factor can be readily transformed from the integer form notation via the one-to-one mapping $i \mapsto e_i$ with domain $\{1, \ldots, K\}$.
We hope this will not cause confusion.
We use the notation $\card x$ to denote the cardinality of some vector $x \in \reals^n$, \ie, the number of nonzero components of $x$.

\subsection{Standard form DLFM fitting problems}\label{sec:dlfm_fit_prob_std}
The model fitting problem of DLFMs can be written as the following mixed integer program:
\begin{equation}\label{prob:std_dlfm}
    \begin{array}{ll}
        \mbox{minimize} & \sum_{i = 1}^{m} z_i^T r_i = \sum_{i = 1}^{m} z_i^T (f(x_i, y_i; \theta_1), \ldots, f(x_i, y_i; \theta_K))\\
        \mbox{subject to} & z_i \in {\{0, 1\}}^K,\quad \card z_i = 1,\quad i = 1, \ldots, m\\
        & \theta_i \in \cvxset,\quad i = 1, \ldots, K,
    \end{array}
\end{equation}
where the problem variables are the latent factors $z_1, \ldots, z_m$ and the model parameters $\theta_1, \ldots, \theta_K$, the problem data are given by feature-observation pairs ${\{x_i, y_i\}}_{i = 1}^m$.
We assume the feasible set $\cvxset$ for model parameters is closed and convex, \ie, can be specified via the DCP ruleset for constraints.
The function $f$ is some metric representing the fitting error, and is assumed to be convex under the DCP ruleset and resolves to a scalar.
The vectors $r_i = (f(x_i, y_i; \theta_1), \ldots, f(x_i, y_i; \theta_K)) \in \reals^K$ are then concatenations of the fitting error under each parameter $\theta_1, \ldots, \theta_K$, evaluated on the dataset ${\{x_i, y_i\}}_{i = 1}^m$.

\paragraph{Extensions.}
Note that the problem (\ref{prob:std_dlfm}) can be readily extended to the cases where the models for generating observations $y$ given feature $x$ are different across latent factors, \ie, consider
\begin{equation*}
    r_i = (f_1(x_i, y_i; \theta_1), \ldots, f_K(x_i, y_i; \theta_K))\quad \mbox{and}\quad \theta_1 \in \cvxset_1,\quad \ldots,\quad \theta_K \in \cvxset_K,
\end{equation*}
for all $i = 1, \ldots, m$.
Theoretically, such an extension does not significantly influence the properties of the problem (\ref{prob:std_dlfm}), as long as our assumptions about the convexity of $f_1, \ldots, f_K$ and $\cvxset_1, \ldots, \cvxset_K$ still hold.
(This can be verified with basic convex analysis.)
Besides, the assumption that the observations are generated via different models under individual latent factors is less considered in practice (since, \eg, there might be some conceptual issues when mixing different metrics together).
Hence, we will focus on the problem formulation (\ref{prob:std_dlfm}) in the subsequent discussion.
In practice, such an extension is indeed technically supported by our implementation in \S\ref{sec:impl}.

\paragraph{Interpretation.}
The DLFM fitting problem (\ref{prob:std_dlfm}) can be interpreted as follows.
For all samples $\{x_i, y_i\}$ in the dataset, given some model parameters $\theta_1, \ldots, \theta_K$, the model fitting error $f(x_i, y_i; \theta_k)$ is evaluated for all $\theta_k$, $k = 1, \ldots, K$, \ie, evaluated under all possible latent factors.
This corresponds to the vector $r_i$, $i = 1, \ldots, m$.
The constraints on the latent factors, given by $z_i \in {\{0, 1\}}^K$ and $\card z_i = 1$ for $i = 1, \ldots, m$, require that the vectors $z_i$ must be one-hot vectors, \ie, in the set of standard basis vectors $\{e_1, \ldots, e_K\}$.
Hence, by forming the inner product $z_i^T r_i$ for all samples in the dataset, we assign a unique label corresponding to a latent factor to each sample, where the model fitting error is cumulated into the objective.
In other word, by solving the problem (\ref{prob:std_dlfm}), we would like to separate the dataset ${\{x_i, y_i\}}_{i = 1}^m$ into $K$ clusters, such that the total model fitting error of all clusters is minimized.

\paragraph{NP-hardness.}
The DLFM fitting problem (\ref{prob:std_dlfm}) is NP-hard (even if we have restricted the function $f$ to be convex).
To show this, consider the $\ell_2$-norm squared metric on $\reals^n$, given by $f(x, y; \theta) = \norm{\theta - x - y}_2^2$.
Let the dataset be ${\{x_i, y_i = 0\}}_{i = 1}^m$, and $\cvxset = \reals^n$.
This leads to an instance of (\ref{prob:std_dlfm}) given by
\begin{equation}\label{prob:dlfm_kmeans}
    \begin{array}{ll}
        \mbox{minimize} & \sum_{i = 1}^{m} z_i^T (\norm{\theta_1 - x_i}_2^2, \ldots, \norm{\theta_K - x_i}_2^2)\\
        \mbox{subject to} & z_i \in {\{0, 1\}}^K,\quad \card z_i = 1,\quad i = 1, \ldots, m
    \end{array}
\end{equation}
with variables $z_1, \ldots, z_m \in \{e_1, \ldots, e_K\}$ and $\theta_1, \ldots, \theta_K \in \reals^n$.
The problem (\ref{prob:dlfm_kmeans}) is equivalent to the $k$-means clustering problem in $\reals^n$ with $K$ clusters.
Hence, the DLFM fitting problem (\ref{prob:std_dlfm}) is at least as hard as the $k$-means clustering problem (\ref{prob:dlfm_kmeans}), which is known to be NP-hard~\cite{aloise2009np}.

\subsection{Examples}
Recall that we assume the function $f$ in (\ref{prob:std_dlfm}) is scalar valued and convex under the DCP ruleset, and the feasible set $\cvxset$ is convex.
Although these assumptions do not support all types of DLFM fitting problems, they still capture a wide range of metric functions for both regression and classification models, including the following examples.

\paragraph{Regression models.}
Consider the function $f$ in (\ref{prob:std_dlfm}) with structure
\begin{equation}\label{eq:loss_reg}
    f(x, y; \theta) = g(x^T \theta - y),
\end{equation}
where $x, \theta \in \reals^n$, $y \in \reals$, and $g \colon \reals \to \reals$ is some loss function, \eg,
\begin{itemize}
    \item square loss: $g(u) = u^2$;
    \item $\ell_p$-loss: $g(u) = \norm{u}_p$ for $p \in [1, \infty]$;
    \item Huber loss: $f(u) = u^2$ for $|u| \leq \delta$, and $f(u) = 2 \delta |u| - \delta^2$ for $|u| > \delta$, where $\delta > 0$ is a parameter.
\end{itemize}
Since $x^T \theta - y$ is affine in $\theta$, the function $f$ given by (\ref{eq:loss_reg}) is convex in $\theta$ if the selected loss function $g$ is convex (which is satisfied by the listed three examples).
The corresponding DLFMs are sometimes named \emph{mixture of linear regressions}.
The same idea can be extended when the observations are in a higher dimensional space.
For example, one may consider the least squares loss $g(u) = \norm{u}_2^2$ or the $\ell_1$-loss $g(u) = \norm{u}_1$ for vector valued observations, and $g(U) = \norm{U}_F^2 = \tr(U^T U)$ (\ie, the square of the Frobenius norm) for matrix valued observations, given that the variable taken by $g$ is affine in $\theta$.

\paragraph{Classification models.}
A simple example for classification DLFMs would be $k$-means clustering, with its corresponding model fitting problem given by (\ref{prob:dlfm_kmeans}).
As a more complex example, consider the loss function $f$ given by
\begin{equation}\label{eq:loss_cls}
    f(X, y; \theta) = -\log\left(\frac{y^T \exp u}{\sum_{i = 1}^{p} \exp u_i}\right),\quad u = X \theta,
\end{equation}
where $X \in \reals^{p \times n}$, $\theta \in \reals^n$, $y \in \{e_1, \ldots, e_p\} \subseteq \reals^p$ is a one-hot label vector.
The intermediate feature vector $u = X \theta$ is then in $\reals^p$ (with its $i$th entry denoted as $u_i$).
The loss function given by (\ref{eq:loss_cls}) is the objective function (\ie, the negative log-likelihood) for multi-class logistic regression problems, which is convex under the DCP ruleset.
To show this, notice that $f$ can be written as
\begin{equation}\label{eq:loss_cls_logsumexp}
    f(X, y; \theta) = -\log\left(\frac{y^T \exp u}{\sum_{i = 1}^{p} \exp u_i}\right) = -\left(y^T u - \log \sum_{i = 1}^{p} \exp u_i \right).
\end{equation}
(To obtain the second equality, we use the fact that $\log (y^T \exp u) = y^T u$ if $y$ is a standard basis vector.)
Since the log-sum-exp expression, given by $\log \sum_{i = 1}^{p} \exp u_i$, is known to be convex in $u$~\cite[\S3.1]{boyd2004convex}, and $u$ is affine in $\theta$, it follows immediately that the function $f$ given by (\ref{eq:loss_cls}) is convex in $\theta$.
The DLFMs corresponding to the problem (\ref{prob:std_dlfm}) with loss function (\ref{eq:loss_cls}) are sometimes referred to as the \emph{hierarchical logistic regression} model.
Note that such formulation includes the binary logistic regression as a special case, and can be readily adapted to deal with hinge loss or exponential loss for binary classification models.

\paragraph{Constraints on parameters.}
In the above examples, we simply assume there is no constraint on the model parameters $\theta$, \ie, $\cvxset = \reals^n$.
It is also common in practice to restrict the parameters, for instance, $\theta \in \reals^n$, to some convex subset $\cvxset$ of $\reals^n$, given some prior information.
The DCP ruleset supports a vast range of constraints that can jointly specify the convex set $\cvxset$, \eg, nonnegativity constraint $\theta \succeq 0$ (where $\succeq$ denotes componentwise inequality), unit norm constraint $\norm{\theta}_2 \leq 1$, and summation constraint $\ones^T \theta = 1$, just list a few.

\section{Heuristic solution method}\label{sec:heuristic}
In this section, we introduce relaxations to the DLFM fitting problem (\ref{prob:std_dlfm}) such that it can be transformed into a multi-convex program and solved (at least approximately).

\subsection{Multi-convex program formulation}
\paragraph{Relaxations.}
We start from relaxing the discrete cardinality constraints about the latent factors $\{z_1, \ldots, z_m\}$ of the problem (\ref{prob:std_dlfm}).
Recall that the constraints, given by
\begin{equation}\label{eq:constr_lf}
    z_i \in {\{0, 1\}}^K,\quad \card z_i = 1,\quad i = 1, \ldots, m,
\end{equation}
restrict the latent factors to be standard basis vectors in $\reals^K$, \ie, $z_i \in \{e_1, \ldots, e_K\}$.
It is in general difficult to handle these constraints when the dataset ${\{x_i, y_i\}}_{i = 1}^m$ is large, and the machinery of combinatorial optimization is often required.
To avoid solving a combinatorial optimization problem, we relax the constraints (\ref{eq:constr_lf}) to
\begin{equation}\label{eq:constr_lf_rel}
    0 \preceq z_i \preceq \ones,\quad \ones^T z_i = 1,\quad i = 1, \ldots, m.
\end{equation}
The relaxed constraints about the latent factor $z_1, \ldots, z_m$, given by (\ref{eq:constr_lf_rel}), are now `smooth' in $z_i$, and can hence be handled by most solving methods.

\paragraph{Interpretation.}
Compared to (\ref{eq:constr_lf}), the constraints given by (\ref{eq:constr_lf_rel}) has the following interpretations:
By solving the DLFM fitting problem with constraints (\ref{eq:constr_lf}), we aim to find a `hard-encoded' (deterministic) label $z_i$ for each sample $\{x_i, y_i\}$ in the dataset, whereas in (\ref{prob:relax_dlfm}), each entry in $z_i$ can be interpreted as the probability of the sample $\{x_i, y_i\}$ belonging to the corresponding factors.
Hence, the relaxed constraints (\ref{eq:constr_lf_rel}) allow `soft-encoded' (stochastic) latent factors for individual samples.
In practice, if required, the stochastic latent factors from (\ref{eq:constr_lf_rel}) can be readily transformed into the deterministic latent factors in (\ref{eq:constr_lf}) via the $\argmax$ operator.

\paragraph{MCP formulation.}
Incorporating the relaxed constraints (\ref{prob:relax_dlfm}), the relaxed DLFM fitting problem is then written as:
\begin{equation}\label{prob:relax_dlfm}
    \begin{array}{ll}
        \mbox{minimize} & \sum_{i = 1}^{m} z_i^T r_i = \sum_{i = 1}^{m} z_i^T (f(x_i, y_i; \theta_1), \ldots, f(x_i, y_i; \theta_K))\\
        \mbox{subject to} & 0 \preceq z_i \preceq \ones,\quad \ones^T z_i = 1,\quad i = 1, \ldots, m\\
        & \theta_i \in \cvxset,\quad i = 1, \ldots, K,
    \end{array}
\end{equation}
with variables $z_1, \ldots, z_m \in \reals^K$ and $\theta_1, \ldots, \theta_K$ depending on the DLFM structure.
We will use the notation $\Theta = (z_1, \ldots, z_m, \theta_1, \ldots, \theta_K)$ to denote the vector of all variables of (\ref{prob:relax_dlfm}) in the sequel.
The relaxed problem (\ref{prob:relax_dlfm}) is multi-convex (or specifically, biconvex) with index sets $\idxset_z$ and $\idxset_\theta$, where $\Theta_{\idxset_z} = (z_1, \ldots, z_m)$ and $\Theta_{\idxset_\theta} = (\theta_1, \ldots, \theta_K)$.
To show this, we need to verify that:
\begin{enumerate}[label=\texttt{(\alph*)}]
    \item\label{cond:sep} $\idxset_z \cap \idxset_\theta = \emptyset$, and
    \item\label{cond:cvx_fixz} the problem (\ref{prob:relax_dlfm}) is convex with $\idxset_z$ fixed, and
    \item\label{cond:cvx_fixtheta} the problem (\ref{prob:relax_dlfm}) is convex with $\idxset_\theta$ fixed.
\end{enumerate}

\paragraph{Proof of multi-convexity.}
Apparently, the condition~\ref{cond:sep} is satisfied since $\Theta_{\idxset_z} \cap \Theta_{\idxset_\theta} = \emptyset$.
To verify the condition~\ref{cond:cvx_fixz}, consider the problem (\ref{prob:relax_dlfm}) with $\idxset_z$ fixed, given by
\begin{equation}\label{prob:relax_dlfm_fixz}
    \begin{array}{ll}
        \mbox{minimize} & \sum_{i = 1}^{m} \tilde{z}_i^T (f(x_i, y_i; \theta_1), \ldots, f(x_i, y_i; \theta_K))\\
        \mbox{subject to} & \theta_i \in \cvxset,\quad i = 1, \ldots, K,
    \end{array}
\end{equation}
with variables $\theta_1, \ldots, \theta_K$ and data ${\{x_i, y_i, \tilde{z}_i\}}_{i = 1}^m$ (where $\tilde{z}_i$ are some fixed latent factors).
Since the loss function $f$ is convex in the variables $\theta_1, \ldots, \theta_K$ and $0 \preceq \tilde{z}_1, \ldots, \tilde{z}_m \preceq \ones$, the objective of (\ref{prob:relax_dlfm_fixz}) is a nonnegative weighted sum of multiple convex functions (and hence convex).
According to our assumption in \S\ref{sec:dlfm_fit_prob_std}, the feasible set $\cvxset$ satisfies the DCP ruleset.
Put together, we conclude that the problem (\ref{prob:relax_dlfm_fixz}) is a convex optimization problem, and thus the condition~\ref{cond:cvx_fixz} is satisfied.
Now we show that the condition~\ref{cond:cvx_fixtheta} also holds.
Similarly, we have the problem (\ref{prob:relax_dlfm}) with $\idxset_\theta$ fixed, given by
\begin{equation}\label{prob:relax_dlfm_fixtheta}
    \begin{array}{ll}
        \mbox{minimize} & \sum_{i = 1}^{m} z_i^T \tilde{r}_i\\
        \mbox{subject to} & 0 \preceq z_i \preceq \ones,\quad \ones^T z_i = 1,\quad i = 1, \ldots, m,
    \end{array}
\end{equation}
with variables $z_1, \ldots, z_m \in \reals^K$.
The data of the problem (\ref{prob:relax_dlfm_fixtheta}) is $\tilde{r}_1, \ldots, \tilde{r}_m \in \reals^K$ obtained from evaluating $f(x_i, y_i; \tilde{\theta}_k)$ for $i = 1, \ldots, m$, $k = 1, \ldots, K$, where $\tilde{\theta}_k$ are some fixed model parameters.
The problem (\ref{prob:relax_dlfm_fixtheta}) is then readily verified to be a linear program, which is, of course, convex.

\subsection{Heuristic solution via block coordinate descent}
A generic method for solving multi-convex optimization problems is the class of \emph{block coordinate descent} (BCD) methods.
The general idea of BCD-type methods (which dates back to Warga~\cite{warga1963minimizing} and Powell~\cite{powell1972search}) is to alternate between multiple convex problems, where each problem is a subproblem of the original multi-convex program with a specific index set fixed.
Specifically, for the relaxed DLFM fitting problem (\ref{prob:relax_dlfm}), in each BCD iteration, we iterate between solving the \emph{parameter (P) problem} and the \emph{factor (F) problem}, given by
\begin{equation}\label{prob:pf}
    \begin{array}{ll}
        \mbox{(P)} &\\
        \mbox{minimize} & \sum_{i = 1}^{m} \tilde{z}_i^T r_i\\
        \mbox{subject to} & r_i = {(f(x_i, y_i; \theta_k))}_{k = 1}^K,\quad \theta_k \in \cvxset\\
        & i = 1, \ldots, m,\quad k = 1, \ldots, K,\\
    \end{array}\qquad
    \begin{array}{ll}
        \mbox{(F)}\\
        \mbox{minimize} & \sum_{i = 1}^{m} z_i^T \tilde{r}_i\\
        \mbox{subject to} & 0 \preceq z_i \preceq \ones,\quad \ones^T z_i = 1\\
        & i = 1, \ldots, m.
    \end{array}
\end{equation}
The P-problem has variables $\theta_1, \ldots, \theta_K$ and data ${\{x_i, y_i\}}_{i = 1}^m$ from the dataset, $\tilde{z}_1, \ldots, \tilde{z}_m \in \reals^K$ corresponding to the optimal point of the F-problem in the last iteration.
Correspondingly, the F-problem has variables $z_1, \ldots, z_m \in \reals^K$ and data $\tilde{r}_i = (f(x_i, y_i; \tilde{\theta}_1), \ldots, f(x_i, y_i; \tilde{\theta}_K))$, $i = 1, \ldots, m$, where $\tilde{\theta}_1, \ldots, \tilde{\theta}_K$ are the optimal point of the P-problem in the last iteration.

\paragraph{Regularizations.}
Note that although it is not mentioned explicitly in (\ref{prob:pf}), regularization terms can be integrated into the objective function of each problem.
As an example, if we would like to obtain some model parameters $\theta_k \in \reals^n$ that are sparse (\ie, have as many zero entries as possible) for all $k = 1, \ldots, K$, we can add an $\ell_1$-regularization term $\lambda \sum_{k = 1}^{K} \norm{\theta_k}_1$ with $\lambda \geq 0$ being the penalty weight to the objective of the P-problem~\cite{tibshirani1996regression}.
As another example regarding the F-problem, if the dataset ${\{x_i, y_i\}}_{i = 1}^m$ are from some time series, a commonly considered assumption for DLFM fitting is that the latent factors should change as less as possible, this leads to the regularization $\lambda \sum_{i = 1}^{m - 1} \Dkl(z_i, z_{i + 1})$ (with weight $\lambda \geq 0$), where for positive vectors $u, v \in \reals_{++}^n$, $\Dkl(u, v) = \sum_{i = 1}^{n} (u_i \log(u_i / v_i) - u_i + v_i)$ is the \emph{Kullback-Leibler (KL) divergence} (which is known to be convex jointly in $(u, v)$~\cite[\S3.2]{boyd2004convex}).
The sparsity of the change points of latent factors comes from the fact that $\sum_{i = 1}^{m - 1} \Dkl(z_i, z_{i + 1}) = \norm{(\Dkl(z_1, z_2), \Dkl(z_2, z_3), \ldots, \Dkl(z_{m - 1}, z_m))}_1$, since the KL-divergence is always nonnegative.

\paragraph{Termination.}
There are several options for determining when to quit the BCD iteration.
The most straightforward and commonly considered approach is to set a maximum iteration number, as used by Shen \etal~\cite{shen2017disciplined}.
Besides, some criteria can be incorporated to automatically terminate the BCD iteration for the two subproblems (\ref{prob:pf}).
For instance, one may quit the algorithm when the gap between the optimal values of the P- and F-problem is less than some very small number $\epsilon \geq 0$ (given that there are no additional regularization terms on both problems), or quit if the update on the variables between iterations is smaller than some threshold.

\paragraph{Convergence.}
In general, very little can be said about the convergence of BCD-type methods for multi-convex problems.
Interested readers may refer to, \eg, Warga~\cite{warga1963minimizing}, Powell~\cite{powell1972search}, and Nutini \etal~\cite{nutini2022let}, for some convergence results under strong assumptions about convexity or differentiability (of the subproblems with some variables fixed).
In case of the relaxed DLFM fitting problem (\ref{prob:relax_dlfm}), one obvious observation is that the objectives of both the P-problem and F-problem are nonincreasing in each BCD iteration, and hence converges.
Despite a lack of strong convergence theory in the general case, the BCD-type methods have been found to be robust and very useful in practice.

\section{Implementation}\label{sec:impl}
In this section, we introduce the implementation of our framework for specifying and solving DLFM fitting problems, based on the \emph{disciplined parameterized programming} (DPP)~\cite{agrawal2019differentiable} implementation of the heuristic solution method discussed in \S\ref{sec:heuristic}.
DPP incorporates symbolic representation to (a subset of) the problem data (\ie, problem \emph{parameters}), such that the value of these parameters can be modified without reconstructing the entire problem.
(Note that the meaning of the word `parameter' for an optimization problem is different from that for a DLFM, as we referred to previously.)
For example, in the case of the problems (\ref{prob:pf}), the parameters are considered to be $\tilde{z}_1, \ldots, \tilde{z}_m \in \reals^K$ for the P-problem, and $\tilde{r}_1, \ldots, \tilde{r}_m \in \reals^K$ for the F-problem.
By specifying a DCP problem according to DPP, solving it repeatedly for different values of the parameters can be much faster than repeatedly solving a new problem.\footnote{For large problems, the non-DPP implementation (which can be readily obtained from our DPP based implementation provided below) is sometimes faster than DPP, since for these problems the canonicalization step in the \texttt{CVXPY} backend may take a very long time.}

The implementation of our framework is given in the listing below.
The corresponding code is fully open-source and is available at
\begin{quote}
    \url{https://github.com/nrgrp/dlfm}.
\end{quote}
\begin{lstlisting}[language=zhpython]
import numpy as np
import cvxpy as cp

### problem data
xs = None  # ndarray: dataset features
ys = None  # ndarray: dataset observations
m = None  # int: number of samples in the dataset

### hyperparameters
eps = 1e-6  # float: termination criterion

### P-problem
K = None  # int: number of latent factors
thetas = []  # list of cp.Variable objects: model parameters
r = []  # list of cp.Expression objects: loss functions

ztil = cp.Parameter((m, K), nonneg=True)
Pobj = cp.sum(cp.multiply(ztil, cp.vstack(r).T))
Preg = 0  # cp.Expression: regularization on model parameters
Pconstr = []  # list of cp.Constraint objects: model parameter constraints
Pprob = cp.Problem(cp.Minimize(Pobj + Preg), Pconstr)
assert Pprob.is_dcp()

### F-problem
rtil = cp.Parameter((K, m))
z = cp.Variable((m, K))
Fobj = cp.sum(cp.multiply(z, rtil.T))
Freg = 0  # cp.Expression: regularization on latent factors
Fconstr = [z >= 0, z <= 1, cp.sum(z, axis=1) == 1]
Fprob = cp.Problem(cp.Minimize(Fobj + Freg), Fconstr)
assert Fprob.is_dcp()

### solve, terminate when the F- and P-objective converge
while True:
    if ztil.value is None:
        ztil.value = np.random.dirichlet(np.ones(K), size=m)
    else:
        ztil.value = np.abs(z.value)
    Pprob.solve()

    rtil.value = cp.vstack(r).value
    Fprob.solve()

    if np.abs(Pobj.value - Fobj.value) < eps:
        break
\end{lstlisting}
To fully specify the fitting problem for customized DLFMs, the users only need to instantiate those objects that are commented in the listing, to which we provide a detailed description in the subsequent paragraphs.

\paragraph{Dataset.}
The dataset for the DLFM fitting problem is specified by the features \texttt{xs}, the observations \texttt{ys}, and the number of samples in the dataset \texttt{m}.
In the general case, the first two objects should be instantiated as a \texttt{numpy.ndarray} with the first dimension equal to the integer \texttt{m}.

\paragraph{P-problem.}
The P-problem requires the user to specify three mandatory objects: \texttt{K}, \texttt{thetas}, \texttt{r}, and two optional objects: \texttt{Preg}, \texttt{Pconstr}.
The lists, \texttt{thetas} and \texttt{r}, represent the parameters and the loss function of the DLFM corresponding to each latent factor, respectively.
The integer \texttt{K} is the number of latent factors.
If regularization and/or constraints need to be applied to the model parameters, the users need to additionally instantiate the objects \texttt{Preg} and/or \texttt{Pconstr}.
The lists \texttt{theta} and \texttt{r} should consist of \texttt{cvxpy.Variable} and \texttt{cvxpy.Expression} objects, respectively.
Note that each expression in the list \texttt{r} has to resolve to a vector in $\reals^K$, according to (\ref{prob:pf}).
In the general case, the integer \texttt{K} should equal to the number of elements in \texttt{thetas} and \texttt{r}.
The regularization term \texttt{Preg} is a \texttt{cvxpy.Expression} object that resolves to a scalar.
The constraints \texttt{Pconstr} is a list of \texttt{cvxpy.Constraint} objects that can be specified according to the standard \texttt{CVXPY} grammar.

\paragraph{F-problem.}
Not much input from the user is required to fully specify the F-problem, except when one would like to integrate regularization terms to the latent factor variables $z_i, \ldots, z_m$ (\ie, the variable \texttt{z} in the code).
In this case, the user needs to instantiate the object \texttt{Freg} as a \texttt{cvxpy.Expression} object that resolves to a scalar (similar to \texttt{Preg} for the P-problem).

\paragraph{Termination.}
We implement the solving iteration as alternating between solving the P- and F-problem, until the gap between the optimal value of the two problems is small enough, controlled by the nonnegative floating point number \texttt{eps}.
The default value of the number \texttt{eps} is set to be $10^{-6}$, which turns out to work well for many DLFM fitting problems (see the numerical examples in \S\ref{sec:example}).
The users can freely change this threshold according to their application scenario.

\section{Examples}\label{sec:example}
This section provides some numerical examples that apply our framework for fitting different DLFMs.
The code corresponding to individual examples is based on the implementation provided in \S\ref{sec:impl} with related lines adapted to respective problems.
Interested readers can access the full implementation at
\begin{quote}
    \url{https://github.com/nrgrp/dlfm}
\end{quote}
as a coding reference.

\subsection{Constrained $k$-means clustering}
\paragraph{Problem description.}
We start from a toy example to demonstrate the basic usage of our framework.
Suppose we are given a dataset consisting of $m = 500$ points with feature dimension $n = 2$, given by $x_1, \ldots, x_m \in \reals^2$, uniformly generated according to $x_i = \bar{x}_i + v_i$, where
\begin{equation}\label{eq:kmeans_dataset}
    \bar{x}_i \in \{x \in \reals^2 \mid \norm{x}_1 = 2\},\quad v_i \sim {\cal N}(0,\ 0.05^2 I),
\end{equation}
and $y_1 = \cdots = y_m = 0$.
(The matrix $I$ is the identity matrix.)
The objective is to partition these points into $K = 4$ clusters such that the within-cluster variance is minimized.
Besides, we require that the center of each cluster $\theta_1, \ldots, \theta_K \in \reals^2$ is constrained in a polyhedron, given by
\begin{equation}\label{eq:kmeans_constr}
    \theta_i \in \{x \in \reals^2 \mid Ax \preceq b\},\quad \mbox{where}\quad A = \left[\begin{array}{rr}
        0.8 & 0.6\\
        -0.7 & 0.9\\
        -1 & -0.5\\
        1 & -1\\
        0.3 & 0.9
    \end{array}\right],\quad
    b = \left[\begin{array}{c}
        1 \\ 0.8 \\ 0.6 \\ 0.7 \\ 0.8
    \end{array}\right].
\end{equation}
This problem corresponds to the $k$-means clustering problem (\ref{prob:dlfm_kmeans}) with linear constraints (\ref{eq:kmeans_constr}).
Formally, the optimization problem of this example is written as
\begin{equation}\label{prob:constr_kmeans}
    \begin{array}{ll}
        \mbox{minimize} & \sum_{i = 1}^{m} z_i^T (\norm{\theta_1 - x_i}_2^2, \ldots, \norm{\theta_K - x_i}_2^2)\\
        \mbox{subject to} & z_i \in {\{0, 1\}}^K,\quad \card z_i = 1,\quad i = 1, \ldots, m\\
        & A\theta_i \preceq b,\quad i = 1, \ldots, K,
    \end{array}
\end{equation}
where $\theta_1, \ldots, \theta_K \in \reals^2$ and $z_1, \ldots, z_m \in \reals^K$ are the variables, $x_1, \ldots, x_m \in \reals^2$, $A \in \reals^{5 \times 2}$ and $b \in \reals^5$ given by (\ref{eq:kmeans_dataset}) and (\ref{eq:kmeans_constr}) are the data.

\paragraph{Problem specification.}
The problems we solve in each BCD iteration, corresponding to the constrained $k$-means clustering problem (\ref{prob:constr_kmeans}), are given by
\begin{equation}\label{prob:kmeans_pf}
    \begin{array}{ll}
        \mbox{(P)} &\\
        \mbox{minimize} & \sum_{i = 1}^{m} \tilde{z}_i^T r_i\\
        \mbox{subject to} & r_i = {\left(\norm{\theta_k - x_i}_2^2\right)}_{k = 1}^K,\quad i = 1, \ldots, m\\
        & A\theta_k \preceq b,\quad k = 1, \ldots, K,\\
    \end{array}\qquad
    \begin{array}{ll}
        \mbox{(F)}\\
        \mbox{minimize} & \sum_{i = 1}^{m} z_i^T \tilde{r}_i\\
        \mbox{subject to} & 0 \preceq z_i \preceq \ones,\quad \ones^T z_i = 1\\
        & i = 1, \ldots, m.
    \end{array}
\end{equation}
The P-problem has variables $\theta_1, \ldots, \theta_K \in \reals^2$ and data $x_1, \ldots, x_m \in \reals^2$, $\tilde{z}_1, \ldots, \tilde{z}_m \in \reals^K$, $A \in \reals^{5 \times 2}$, $b \in \reals^5$.
The F-problem has variables $z_1, \ldots, z_m \in \reals^K$ and data $\tilde{r}_1, \ldots, \tilde{r}_m \in \reals^K$.
To type the problems (\ref{prob:kmeans_pf}) into our framework, the major adaptations we need to make are given by the following code.
\begin{verbatim}
    for k in range(K):
        thetas.append(cp.Variable((1, n)))
        r.append(cp.sum(cp.square(xs - thetas[-1]), axis=1))
    Pconstr = [A @ theta.T <= b for theta in thetas]
\end{verbatim}

\paragraph{Numerical result.}
\begin{figure}[t]
    \centering
    \includegraphics[width=0.7\textwidth]{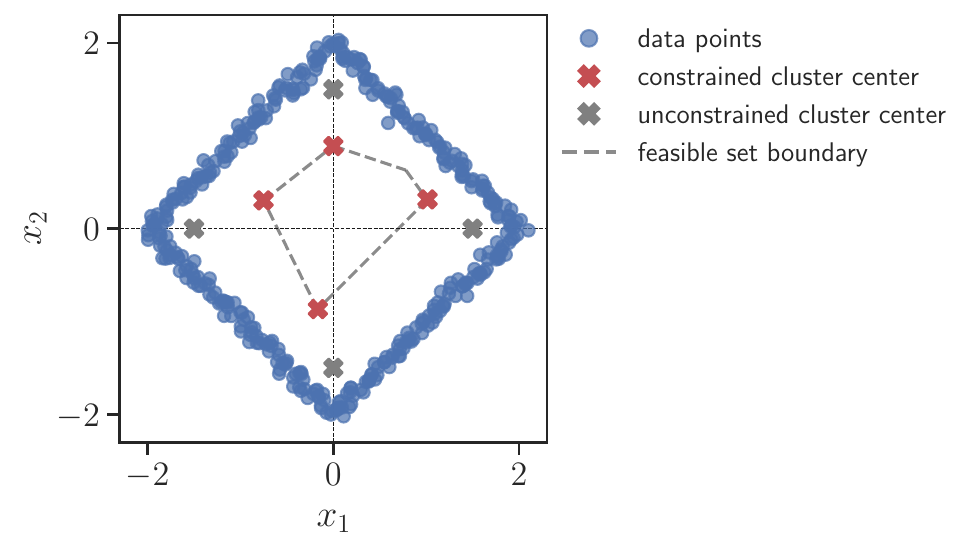}
    \caption{
        Results for the constrained $k$-means clustering example.
    }\label{fig:constr_kmeans}
\end{figure}

The result for this example is shown in figure~\ref{fig:constr_kmeans}.
The blue dots represent the dataset, and the gray dashed lines are the boundary of the feasible set given by (\ref{eq:kmeans_constr}).
The returned centers for different clusters are shown as red crosses, which are all at the vertices of the polyhedron defined by the feasible set.
As a reference, the cluster centers for the same dataset when removing the constraints (\ref{eq:kmeans_constr}) are plotted in gray crosses.

\subsection{Mixture of linear regressions}
\paragraph{Problem description.}
We now demonstrate an application of our framework in fitting a mixture of linear regressions.
Such a DLFM model has a long application history in machine learning~\cite{bengio1994input,gaffney1999trajectory,li2018learning}.
Suppose we are given a dataset ${\{x_i, y_i\}}_{i = 1}^m$ generated according to
\begin{equation*}
    x_i \sim {\cal U}(a, b),\quad
    \theta_i \sim {\rm Cat}(\{\theta_1, \ldots, \theta_K\},\ p),\quad
    y_i \sim {\cal N}(x_i^T \theta_i, \sigma^2),
\end{equation*}
\ie, for the $i$th sample, the feature vector $x_i \in \reals^n$ is first generated from a uniform distribution between the lower bound $a \in \reals^n$ and upper bound $b \in \reals^n$ (assume $a \prec b$), then the linear combination coefficients $\theta_i \in \reals^n$ are generated from a categorical distribution on the set $\{\theta_1, \ldots, \theta_K\}$ according to probabilities $p \in \reals_+^K$ with $\ones^T p = 1$, and finally the observation $y_i \in \reals$ is generated from the Gaussian distribution with mean $x_i^T \theta_i$ and variance $\sigma^2$.
To fit this mixture of linear regressions given the dataset ${\{x_i, y_i\}}_{i = 1}^m$, we need to recover the set of parameters $\{\theta_1, \ldots, \theta_K\}$, as well as the latent factor labels indicating which $\theta_i \in \{\theta_1, \ldots, \theta_K\}$ was used to generate the observation $y_i$ given feature $x_i$, for all $i = 1, \ldots, m$.
This optimization problem in standard form (\ref{prob:std_dlfm}) is given by
\begin{equation}\label{prob:mix_lin_reg}
    \begin{array}{ll}
        \mbox{minimize} & \sum_{i = 1}^{m} z_i^T ({(x_i^T \theta_1 - y_i)}^2, \ldots, {(x_i^T \theta_K - y_i)}^2)\\
        \mbox{subject to} & z_i \in {\{0, 1\}}^K,\quad \card z_i = 1,\quad i = 1, \ldots, m,
    \end{array}
\end{equation}
where $\theta_1, \ldots, \theta_K \in \reals^n$ and $z_1, \ldots, z_m \in \reals^K$ are the variables, $x_1, \ldots, x_m \in \reals^n$ and $y_1, \ldots, y_m \in \reals$ are the data.

\paragraph{Problem specification.}
The problems we solve in each BCD iteration corresponding to (\ref{prob:mix_lin_reg}) are given by
\begin{equation}\label{prob:mix_lin_reg_pf}
    \begin{array}{ll}
        \mbox{(P)} &\\
        \mbox{minimize} & \sum_{i = 1}^{m} \tilde{z}_i^T r_i\\
        \mbox{subject to} & r_i = {\left({(x_i^T \theta_k - y_i)}^2\right)}_{k = 1}^K\\
        & i = 1, \ldots, m,
    \end{array}\qquad
    \begin{array}{ll}
        \mbox{(F)}\\
        \mbox{minimize} & \sum_{i = 1}^{m} z_i^T \tilde{r}_i\\
        \mbox{subject to} & 0 \preceq z_i \preceq \ones,\quad \ones^T z_i = 1\\
        & i = 1, \ldots, m.
    \end{array}
\end{equation}
The P-problem has variables $\theta_1, \ldots, \theta_K \in \reals^n$ and data $x_1, \ldots, x_m \in \reals^n$, $y_1, \ldots, y_m \in \reals$, $\tilde{z}_1, \ldots, \tilde{z}_m \in \reals^K$; the F-problem has variables $z_1, \ldots, z_m \in \reals^K$ and data $\tilde{r}_1, \ldots, \tilde{r}_m \in \reals^K$.
The major adaptations we need to make to specify the problems (\ref{prob:mix_lin_reg_pf}) using our framework are given by the following code.
\begin{verbatim}
    for k in range(K):
        thetas.append(cp.Variable(n))
        r.append(cp.square(xs @ thetas[-1] - ys))
\end{verbatim}

\paragraph{Numerical result.}
\begin{figure}[t]
    \centering
    \includegraphics[width=0.4\textwidth]{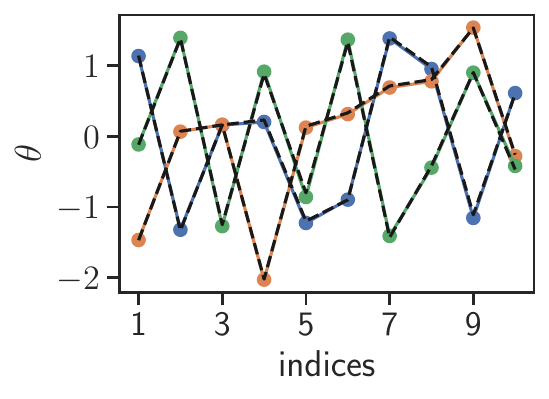}
    \caption{
        Results for the mixture of linear regressions example.
        The three colored solid lines are the recovered parameters and the black dashed lines are the corresponding ground truth.
    }\label{fig:mix_lin_reg}
\end{figure}

The dataset put into test in this example consists of $m = 500$ samples in $\reals^{10}$ with $K = 3$, and was generated under the following parameters:
\begin{align*}
    & a = (-10,\ \ldots,\ -10) \in \reals^{10},\quad b = (10,\ \ldots,\ 10) \in \reals^{10}\\
    & \theta_1 = (-1.47,\ 0.07,\ 0.16,\ -2.02,\ 0.14,\ 0.33,\ 0.71,\ 0.80,\ 1.53,\ -0.26)\\
    & \theta_2 = (-0.12,\ 1.38,\ -1.25,\ 0.88,\ -0.80,\ 1.33,\ -1.43,\ -0.42,\ 0.90,\ -0.47)\\
    & \theta_3 = (1.14,\ -1.33,\ 0.16,\ 0.23,\ -1.20,\ -0.90,\ 1.40,\ 0.98,\ -1.11,\ 0.60)\\
    & p = (0.4,\ 0.3,\ 0.3),\quad \sigma^2 = 1.5^2.
\end{align*}
The recovered set of parameters is shown in figure~\ref{fig:mix_lin_reg}, where the solid lines are the fitted parameters and the black dashed lines are the ground truth.
It can be seen that the recovered and true parameters align exactly.
Besides, the accuracy of recovering the latent factor labels for each sample in the dataset is $0.94$.

\subsection{Hierarchical forgetting Q-learning}
\paragraph{Problem description.}
We consider the \emph{hierarchical forgetting Q-learning} model that is widely used in neuroscience and psychology for the modeling of subject's decision making behavior under multi-armed bandit tasks~\cite{ito2009validation,beron2022mice,zhu2024multiintention,zhu2025solving}.
Consider an agent repeatedly facing a choice among $p$ actions, with (potentially different) probabilities of obtaining a reward from each action.
After the choice at time step $t - 1$, the agent receives a \emph{reward signal} $u(t) \in \reals^p$ that depends on the selected action, \eg,
\begin{equation}\label{eq:rew_signal}
    u_i(t) = \left\{
        \begin{array}{ll}
            1 & \mbox{if action $i$ was selected \emph{and} rewarded}\\
            0 & \mbox{otherwise},
        \end{array}\right.
\end{equation}
for $i = 1, \ldots, p$, where $u_i(t)$ denotes the $i$th element of the vector $u(t)$.
To select the action for the next time step $t$, the agent formulates a \emph{value function} $v(t) \in \reals^p$ according to
\begin{equation}
    v(t) = X(t) \theta(t),\quad
    X(t) = \left[
        \begin{array}{cccc}
            u(t) & u(t - 1) & \cdots & u(t - n + 1)
        \end{array}
    \right] \in \reals^{p \times n},
\end{equation}
where $\theta(t) \in \{\theta_1, \ldots, \theta_K\} \subseteq \reals^n$ is the linear combination coefficient, and the vectors $u(t - \tau + 1)$ with $t - \tau + 1 \leq 0$ for $\tau = 1, \ldots, n$ are padded with zero.
Then the action $y(t) \in \{e_1, \ldots, e_p\} \subseteq \reals^p$ is selected according to
\begin{equation}\label{eq:pi}
    y(t) \sim {\rm Cat}(\{e_1, \ldots, e_p\},\ \exp v(t) / \ones^T \exp v(t)).
\end{equation}
The fitting problem of a hierarchical forgetting Q-learning model consists in recovering the parameters $\{\theta_1, \ldots, \theta_K\}$ as well as the latent factor labels indicating which $\theta(t) \in \{\theta_1, \ldots, \theta_K\}$ was used at time step $t$, for all $t = 1, \ldots, m$, given a dataset ${\{X(t), y(t)\}}_{t = 1}^m$ observed from an agent following (\ref{eq:rew_signal}) to (\ref{eq:pi}).

\paragraph{Dataset and problem formulation.}
In this example, we consider a multi-armed bandit environment with the number of arms $p = 3$ and reward probabilities $(0.1,\ 0.2,\ 0.7)$.
The agent selects its action at each time step following (\ref{eq:rew_signal}) to (\ref{eq:pi}) with $\theta(t) \in \{\theta_1, \theta_2\} \subseteq \reals^5$ (\ie, $K = 2$, $n = 5$), given by
\begin{equation*}
    \begin{array}{rcl}
        \theta_1 & = & (9.9,\ 9.9 \times 10^{-2},\ 9.9 \times 10^{-4},\ 9.9 \times 10^{-6},\ 9.9 \times 10^{-8})\\
        \theta_2 & = & (-4,\ -0.8,\ -0.16,\ -0.032,\ -0.0064).
    \end{array}
\end{equation*}
These two groups of parameters can be interpreted as follows.
Under parameters $\theta_1$, the agent tends to exploit the action that was mostly rewarded in the last $5$ trials, whereas under parameters $\theta_2$, the agent tends to explore the action that had the least reward in the last $5$ trials.
In the fitting problem, suppose we are given the prior information that $\theta_1$ is nonnegative and nonincreasing, and $\theta_2$ is nonpositive and nondecreasing.
Such prior information corresponds to the following constraints:
\begin{equation*}
        \theta_1 \geq 0,\quad \theta_{1, 1} \geq \cdots \geq \theta_{1, 5},\quad \theta_2 \leq 0,\quad \theta_{2, 1} \leq \cdots \leq \theta_{2, 5},
\end{equation*}
where $\theta_{1, i}$ denotes the $i$th entry of the vector $\theta_1$.
The dataset consists of $m = 200$ trials, where the agent starts with $\theta_1$ and switches to another combination coefficient every $20$ trials.
Such a latent factor transition dynamics indicates that a regularization term on the latent factor labels should be applied such that its transition is sparse.
Put together, the optimization problem for this example is given by
\begin{equation}\label{prob:hier_forget_q_learning}
    \begin{array}{ll}
        \mbox{minimize} & -\sum_{t = 1}^{m} {z(t)}^T \log \left(\frac{{y(t)}^T \exp(X(t)\theta_1)}{\ones^T \exp(X(t)\theta_1)},\ \frac{{y(t)}^T \exp(X(t)\theta_2)}{\ones^T \exp(X(t)\theta_2)}\right)\\
        &\qquad + \lambda \sum_{t = 1}^{m - 1} \Dkl(z(t), z(t + 1))\\
        \mbox{subject to} & z(t) \in {\{0, 1\}}^2,\quad \card z(t) = 1,\quad i = 1, \ldots, m\\
        & \theta_1 \geq 0,\quad \theta_{1, 1} \geq \cdots \geq \theta_{1, 5}\\
        &  \theta_2 \leq 0,\quad \theta_{2, 1} \leq \cdots \leq \theta_{2, 5}
    \end{array}
\end{equation}
(\cf, (\ref{eq:loss_cls}) and (\ref{eq:loss_cls_logsumexp})).
The problem (\ref{prob:hier_forget_q_learning}) has variables $\theta_1, \theta_2 \in \reals^5$, $z(1), \ldots, z(m) \in \reals^2$, and data $X(1), \ldots, X(m) \in \reals^{3 \times 5}$, $y(1), \ldots, y(m) \in \reals^3$.
The regularization weight $\lambda \geq 0$ is the hyperparameter.

\paragraph{Problem specification.}
The problems we solve in each BCD iteration corresponding to (\ref{prob:hier_forget_q_learning}) are given by
\begin{equation}\label{prob:hier_forget_q_learning_pf}
    \begin{aligned}
        &\mbox{(P)}\quad
        \begin{array}{ll}
            \mbox{minimize} & \sum_{t = 1}^{m} {\tilde{z}(t)}^T r(t)\\
            \mbox{subject to} & r(t) = -\log \left(\frac{{y(t)}^T \exp(X(t)\theta_1)}{\ones^T \exp(X(t)\theta_1)},\ \frac{{y(t)}^T \exp(X(t)\theta_2)}{\ones^T \exp(X(t)\theta_2)}\right),\quad t = 1, \ldots, m\\
            & \theta_1 \geq 0,\quad \theta_{1, 1} \geq \cdots \geq \theta_{1, 5}\\
            & \theta_2 \leq 0,\quad \theta_{2, 1} \leq \cdots \leq \theta_{2, 5},
        \end{array}\\[5pt]
        &\mbox{(F)}\quad
        \begin{array}{ll}
            \mbox{minimize} & \sum_{t = 1}^{m} {z(t)}^T \tilde{r}(t) + \lambda \sum_{t = 1}^{m - 1} \Dkl(z(t), z(t + 1))\\
            \mbox{subject to} & 0 \preceq z(t) \preceq \ones,\quad \ones^T z(t) = 1,\quad t = 1, \ldots, m.
        \end{array}
    \end{aligned}
\end{equation}
The P-problem has variables $\theta_1, \theta_2 \in \reals^5$ and data $X(1), \ldots, X(m) \in \reals^{3 \times 5}$, $y(1), \ldots, y(m) \in \reals^3$, $\tilde{z}(1), \ldots, \tilde{z}(m) \in \reals^2$.
The F-problem has variables $z(1), \ldots, z(m) \in \reals^2$ and data $\tilde{r}(1), \ldots, \tilde{r}(m) \in \reals^2$.
The major adaptations we need to make to specify the problems (\ref{prob:hier_forget_q_learning_pf}) using our framework are given by the following code.
\begin{verbatim}
    for k in range(K):
        thetas.append(cp.Variable(n))
        r.append(cp.hstack([-(xs[i] @ thetas[-1] @ ys[i] -
                    cp.log_sum_exp(xs[i] @ thetas[-1])) for i in range(m)]))
    Pconstr = [thetas[0] >= 0, cp.diff(thetas[0]) <= 0,
               thetas[1] <= 0, cp.diff(thetas[1]) >= 0]
    Freg = lbd * cp.sum(cp.kl_div(z[:-1], z[1:]))
\end{verbatim}

\paragraph{Numerical result.}
\begin{figure}[t]
    \centering
    \includegraphics[width=0.9\textwidth]{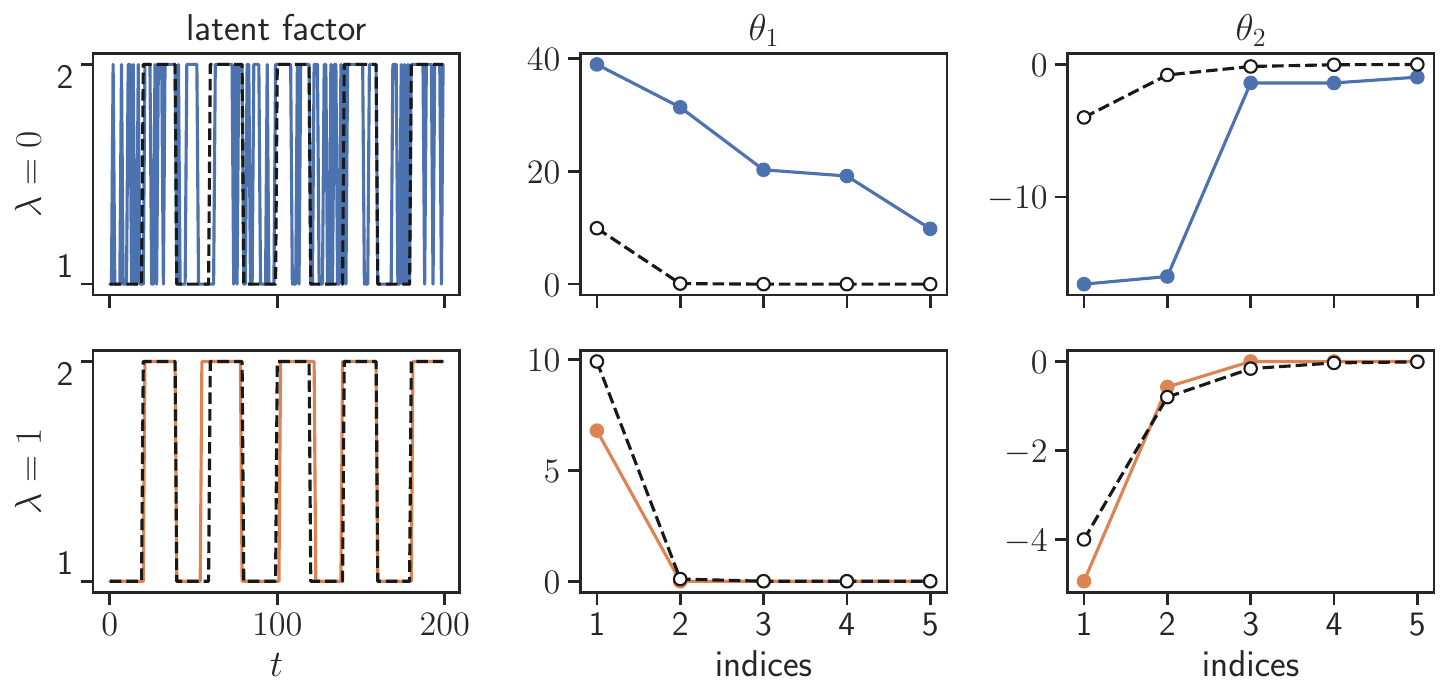}
    \caption{
        Results for the hierarchical forgetting Q-learning example under different $\lambda$ values.
        The colored solid lines are the recovered latent factor labels (left) and model parameters (middle and right).
        The black dashed lines are the corresponding ground truth.
    }\label{fig:hier_forget_q_learning}
\end{figure}

We tested the fitting performance under two different regularization weights, \ie, $\lambda = 0$ and $\lambda = 1$, where the former condition is equivalent to removing the smoothness regularization term on $z(t)$.
The results are shown in figure~\ref{fig:hier_forget_q_learning}.
It can be seen directly that if no regularization terms are added ($\lambda = 0$), the recovered latent factor labels switch drastically (figure~\ref{fig:hier_forget_q_learning}, top left), which leads to only a $0.71$ accuracy of recovering the correct label.
As a result, the recovered model parameters $\theta_1$ and $\theta_2$ are far from the ground truth (figure~\ref{fig:hier_forget_q_learning}, top middle and right).
When we solve the problems (\ref{prob:hier_forget_q_learning_pf}) with hyperparameter $\lambda = 1$, on the contrary, the recovered latent factor labels align almost exactly to the ground truth (figure~\ref{fig:hier_forget_q_learning}, bottom left), resulting in an $0.93$ accuracy of latent factor label identification.
Correspondingly, the recovered model parameters $\theta_1$ and $\theta_2$ are much closer to the ground truth value (figure~\ref{fig:hier_forget_q_learning}, bottom middle and right).

\subsection{Input-output hidden Markov model}
\paragraph{Problem description.}
We consider a special case of the \emph{input-output hidden Markov model} (IO-HMM)~\cite{bengio1994input} with linear outputs, which has been widely used in neuroscience for behavior modeling and neural activities analysis~\cite{escola2011hidden,calhoun2019unsupervised,bolkan2022opponent,ashwood2022mice,jha2024active}.
Let $\hat{z} \in \{1, \ldots, K\}$ be the latent factor label of a $K$-state IO-HMM with initial state distribution $p_{\rm init} \in \reals^K$ with $\ones^T p_{\rm init} = 1$ and transition matrix $P_{\rm tr} \in \reals^{K \times K}$ with $P_{\rm tr} \ones = \ones$.
At the time step $t$, the latent factor label $\hat{z}(t)$ is sampled according to
\begin{equation*}
    \hat{z}(t) \sim \left\{
        \begin{array}{ll}
            {\rm Cat}(p_{\rm init}) & t = 0\\
            {\rm Cat}(p_{\hat{z}(t - 1)}) & t > 0,
        \end{array}\right.
\end{equation*}
where the vector $p_{\hat{z}(t-1)} \in \reals^K$ denotes the $\hat{z}(t-1)$th row of the matrix $P_{\rm tr}$.
The input feature vector $x(t) \in \reals^n$ to the IO-HMM at this time step is generated according to
\begin{equation}\label{eq:io_hmm_feature}
    x(t) = (\bar{x}(t),\ 1),\quad \bar{x}(t) \sim {\cal U}(a, b),
\end{equation}
where $a, b \in \reals^{n - 1}$ ($a \prec b$) are the lower and upper bounds of the uniform distribution.
Note that in (\ref{eq:io_hmm_feature}), we implicitly integrate a bias term into the last entry of each $x(t)$.
The output $y(t) \in \{0, 1\}$ of this IO-HMM at time step $t$ is then generated from a logistic model, \ie,
\begin{equation*}
    \prob(y(t) = 1) = \frac{1}{1 + \exp(-{x(t)}^T \theta_{\hat{z}(t)})},
\end{equation*}
where $\theta_{\hat{z}(t)} \in \{\theta_1, \ldots, \theta_K\} \subseteq \reals^n$ is the linear combination coefficient.
For the fitting problem of this IO-HMM, the objective is to recover the Markov process transition matrix $P_{\rm tr}$ (and the initial state distribution $p_{\rm init}$ if multiple observed sequences are given), the model parameters $\theta_1, \ldots, \theta_K$, and the unobserved latent factor labels $\hat{z}(1), \ldots, \hat{z}(m)$, given the dataset ${\{x(t), y(t)\}}_{t = 1}^m$.

\paragraph{Dataset and problem formulation.}
In this example, we consider an IO-HMM with $K = 3$ and $n = 2$, given by the following parameters:
\begin{align*}
    &a = -5,\quad b = 5\\
    &\theta_1 = (-2, 0),\quad \theta_2 = (2, 6),\quad \theta_3 = (3, -5)\\
    &p_{\rm init} = (1, 0, 0),\quad 
    P_{\rm tr} = \left[\begin{array}{lll}
        0.90 & 0.05 & 0.05\\
        0.01 & 0.98 & 0.01\\
        0.03 & 0.02 & 0.95
    \end{array}\right].
\end{align*}
The dataset ${\{x(t), y(t)\}}_{t = 1}^m$ is a single sequence observed from this IO-HMM with $m = 500$ samples.
To formulate the optimization problem, we assume that prior information about the model parameters $\theta_1, \theta_2, \theta_3$ is given by the following constraints:
\begin{equation*}
    \theta_{1, 1} \leq 0,\quad \theta_{2, 1} \geq 0,\quad \theta_{3, 1} \geq 0.
\end{equation*}
Besides, we also integrate an $\ell_2$-regularization term on the model parameters, and a smoothness regularization term on the latent factor labels as in (\ref{prob:hier_forget_q_learning}).
Put together, we have
\begin{equation}\label{prob:io_hmm}
    \begin{array}{ll}
        \mbox{minimize} & -\sum_{t = 1}^{m} {z(t)}^T {\left(y(t){x(t)}^T \theta_k - \log(1 + \exp({x(t)}^T \theta_k))\right)}_{k = 1}^{3}\\
        &\qquad + \lambda_\theta \sum_{k = 1}^{3} \norm{\theta_k}_2 + \lambda_z \sum_{t = 1}^{m - 1} \Dkl(z(t), z(t + 1))\\
        \mbox{subject to} & z(t) \in {\{0, 1\}}^3,\quad \card z(t) = 1,\quad i = 1, \ldots, m\\
        & \theta_{1, 1} \leq 0,\quad \theta_{2, 1} \geq 0,\quad \theta_{3, 1} \geq 0,
    \end{array}
\end{equation}
where $\theta_1, \theta_2, \theta_3 \in \reals^2$ and $z(1), \ldots, z(m) \in \reals^3$ are the problem variables, $x(1), \ldots, x(m) \in \reals^2$ and $y(1), \ldots, y(m) \in \{0, 1\}$ are the problem data, $\lambda_\theta \geq 0$ and $\lambda_z \geq 0$ are the hyperparameters.
Note that the objective of recovering the transition matrix $P_{\rm tr}$ is not included in (\ref{prob:io_hmm}).
Instead, the transition matrix $P_{\rm tr}$ is estimated according to the returned latent factor labels $z(1), \ldots, z(m)$ after solving (\ref{prob:io_hmm}).

\paragraph{Problem specification.}
The problems we solve in each BCD iteration corresponding to (\ref{prob:io_hmm}) are given by
\begin{equation}\label{prob:io_hmm_pf}
    \begin{aligned}
        &\mbox{(P)}\quad
        \begin{array}{ll}
            \mbox{minimize} & \sum_{t = 1}^{m} {\tilde{z}(t)}^T r(t) + \lambda_\theta \sum_{k = 1}^{3} \norm{\theta_k}_2\\
            \mbox{subject to} & r(t) = -{\left(y(t){x(t)}^T \theta_k - \log(1 + \exp({x(t)}^T \theta_k))\right)}_{k = 1}^{3},\quad t = 1, \ldots, m\\
            & \theta_{1, 1} \leq 0,\quad \theta_{2, 1} \geq 0,\quad \theta_{3, 1} \geq 0,
        \end{array}\\[5pt]
        &\mbox{(F)}\quad
        \begin{array}{ll}
            \mbox{minimize} & \sum_{t = 1}^{m} {z(t)}^T \tilde{r}(t) + \lambda_z \sum_{t = 1}^{m - 1} \Dkl(z(t), z(t + 1))\\
            \mbox{subject to} & 0 \preceq z(t) \preceq \ones,\quad \ones^T z(t) = 1,\quad t = 1, \ldots, m.
        \end{array}
    \end{aligned}
\end{equation}
The P-problem has variables $\theta_1, \theta_2, \theta_3 \in \reals^2$ and data $x(1), \ldots, x(m) \in \reals^2$, $y(1), \ldots, y(m) \in \{0, 1\}$, $\tilde{z}(1), \ldots, \tilde{z}(m) \in \reals^3$.
The F-problem has variables $z(1), \ldots, z(m) \in \reals^3$ and data $\tilde{r}(1), \ldots, \tilde{r}(m) \in \reals^3$.
The major adaptations we need to make to specify the problems (\ref{prob:hier_forget_q_learning_pf}) using our framework are given by the following code.
\begin{verbatim}
    for k in range(K):
        thetas.append(cp.Variable(n))
        r.append(-(cp.multiply(ys, xs @ thetas[-1]) - 
                    cp.logistic(xs @ thetas[-1])))
    Preg = lbd_theta * cp.sum(cp.norm2(cp.vstack(thetas), axis=1))
    Pconstr = [thetas[0][0] <= 0, thetas[1][0] >= 0, thetas[2][0] >= 0]
    Freg = lbd_z * cp.sum(cp.kl_div(z[:-1], z[1:]))
\end{verbatim}

\paragraph{Numerical result.}
\begin{figure}
    \centering
    \includegraphics[width=0.75\textwidth]{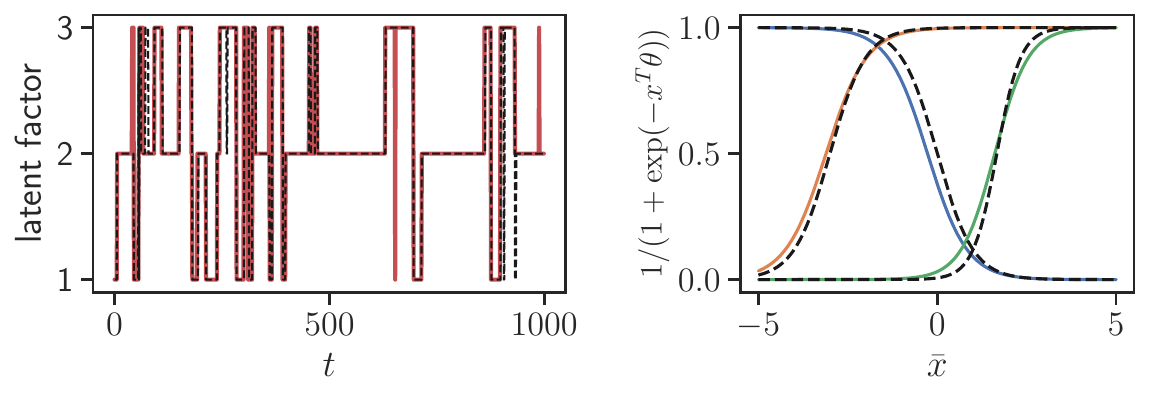}
    \caption{
        Results for the input-output hidden Markov model example.
        The colored solid lines represent the recovered latent factor labels (left) and the decision curve under parameters $\theta_1, \theta_2, \theta_3$, respectively (right).
        The black dashed lines are the corresponding ground truth.
    }\label{fig:io_hmm}
\end{figure}

We specified and solved the problems (\ref{prob:io_hmm_pf}) with hyperparameters $\lambda_\theta = 0.5$ and $\lambda_z = 1$.
The results are shown in figure~\ref{fig:io_hmm}.
The recovered transition matrix $P_{\rm tr}$ is
\begin{equation*}
    \left[\begin{array}{lll}
        0.920 & 0.044 & 0.036\\
        0.003 & 0.978 & 0.019\\
        0.030 & 0.030 & 0.940
    \end{array}\right].
\end{equation*}

\section*{Acknowledgments}
This work has been funded as part of BrainLinks-BrainTools, which is funded by the Federal Ministry of Economics, Science and Arts of Baden-W\"urttemberg within the sustainability program for projects of the Excellence Initiative II, and CRC/TRR 384 ``IN-CODE''.

\newpage
\bibliography{refs}

\end{document}